\documentclass[12pt]{amsart}
\usepackage{amsfonts, amsmath, amsthm}
\usepackage{epsfig}
\usepackage{graphicx}

\providecommand{\bysame}{\leavevmode\hbox to3em{\hrulefill}\thinspace}
\providecommand{\MR}{\relax\ifhmode\unskip\space\fi MR }
\providecommand{\MRhref}[2]{%
  \href{http://www.ams.org/mathscinet-getitem?mr=#1}{#2}
}
\providecommand{\href}[2]{#2}

\newtheorem{pro}{Proposition}[section]
\newtheorem{thm}[pro]{Theorem}
\newtheorem{lem}[pro]{Lemma}

\newtheorem{question}[pro]{Question}

\newtheorem{cnj}[pro]{Conjecture}

\newtheorem{rmkks}[pro]{Remarks}

\newtheorem{cor}[pro]{Corollary}

\theoremstyle{definition}

\theoremstyle{remark}
\newtheorem*{rmk}{Remark}

\newcommand{\scc}{simple closed curve}
\newcommand{\s}{\Sigma}

\newcommand{\ie}{{\it i.e.}}

\newcommand{\del}{\partial}

\title{On the growth rate of tunnel number of knots}

\date{\today}
\address{Department of Mathematics, Nara Women's University
Kitauoya Nishimachi, Nara 630-8506, Japan}
\address{Department of mathematical Sciences, University of
Arkansas, Fayetteville, AR 72701}
\email{tsuyoshi@cc.nara-wu.ac.jp} \email{yoav@uark.edu}
\author{Tsuyoshi Kobayashi}
\author{Yo'av Rieck}
\thanks{The first named author was supported in part by grant in aid JSPS 15540073}

\begin{document}

\maketitle

\begin{abstract}

Given a knot $K$ in a closed orientable manifold $M$ we define the
growth rate of the tunnel number of $K$ to be $gr_t(K) =
\limsup_{n \to \infty} \frac{t(nK) - n t(K)}{n-1}$. As our main
result we prove that the Heegaard genus of $M$ is strictly less
than the Heegaard genus of the knot exterior if and only if the
growth rate is less than 1. In particular this shows that a
non-trivial knot in $S^3$ is never asymptotically super additive.
The main result gives conditions that imply falsehood of
Morimoto's Conjecture.
\end{abstract}

\section{Introduction}
\label{sec:intro}

This paper is concerned with the asymptotic behavior of the tunnel
number of knots in closed 3-manifolds under repeated connected sum
operation.  Let $t(K)$ be the tunnel number of a knot $K$ (for
definitions see Section \ref{sec:prelims}).  It is well known that
for any pair of knots $K_1$, $K_2$ the following inequality holds:

\begin{equation}
\label{eq:trivial} t(K_1 \# K_2) \leq t(K_1) + t(K_2) +1.
\end{equation}

For convenience we denote $\#_{i=1}^n K$, the connected sum  of
$n$ copies of $K$, by $nK$.  By applying Inequality
(\ref{eq:trivial}) repeatedly we obtain:

\begin{equation}
\label{eq:trivials} t(n K) \leq n t(K) + (n-1).
\end{equation}

We define the \em growth rate \em of the tunnel number of $K$,
denoted by $gr_t(K)$, to be:

\begin{equation}
\label{eq:growth}
  gr_t(K) = \limsup_{m \to \infty}  \frac{t(mK) - mt(K)}{m-1}.
\end{equation}

\begin{rmkks}{\rm
\label{rmks:basic facts}

\begin{enumerate}
    \item Inequality (\ref{eq:trivial}) implies that $\frac{t(mK) -
    mt(K)}{m-1} \le 1$ for any $m$.  Since $t(mK) \geq 0$, we obtain
    $\frac{t(mK) - mt(K)}{m-1} \geq \frac{- m}{m-1}t(K)$.  Combining
    the two we see that $-t(K) \leq gr_t(K) \leq 1$. Therefore
    $gr_t(K)$ exists and is finite.
    \item We say that a knot $K$ is \em meridionally small \em if
    the meridian of $K$ is not the boundary of an essential surface (\ie,
    $\pi_1$ injective, not boundary parallel surface) embedded in $E(K)$.
    If $K$ is meridionally small then:

    $$gr_t(K) \geq 0.$$

    In fact, if $K$ is meridionally small then, by
    \cite{morimoto2}and \cite{kobayashi-rieck}, we have $t(nK)
    \geq nt(K)$ for any $n$, which implies the above inequality.
    \item If $K$ is a knot in $S^3$ then:

    $$gr_t(K) \geq -\frac{2}{3}t(K) - 1.$$

    In fact, Scharlemann--Schultens \cite{schsch} (quoting Kwong
    \cite{kwong}) have the inequality $t(nK) \geq \frac{n}{3}t(K)
    - (n-1)$.  This implies:

    \begin{eqnarray*}
      \frac{t(nK) - nt(K)}{n-1} & \geq &  \frac{(n/3) t(K) - (n-1) - nt(K)}{n-1}\\
       &=& \frac{(-2n/3) t(K)-(n-1)}{n-1}.\\
    \end{eqnarray*}

    Taking limit gives the above inequality.

\end{enumerate}

}
\end{rmkks}

Note that for any knot $K$, the Heegaard genus of the knot
exterior $E(K)$, denoted $g(E(K))$, is equal to $t(K)+1$, and it
is at least the Heegaard genus of the ambient manifold $M$, denote
$g(M)$. In this paper, we prove the following (for definitions see
Section \ref{sec:prelims}):

\begin{thm}
\label{thm:main}

Let $K$ be a knot in a closed, orientable 3-manifold $M$ such that
the $g(M) < g(E(K))$. Suppose $K$ is a genus $t(K)$, $n$-bridge
knot (\ie, $n$ is the bridge index of $K$ with respect to Heegaard
surfaces of genus $t(K)$ $(=g(E(K)) - 1)$). Then $gr_t(K) \le
\frac{n-1}{n}$. In particular, for sufficiently large $n$ equality
does not hold in Inequality (\ref{eq:trivials}).
\end{thm}

Note that the assumption of Theorem \ref{thm:main} holds trivially
for any knot in $S^3$.

Suppose that a knot $K$ is isotopic onto a genus $g$ Heegaard
surface $\s \subset M$.  (We are not assuming that $\s$ is of
minimal genus.)  By connecting a spine of one of the handlebodies
complementary to $\s$ to $K$ via a vertical arc it is easy to see
that $t(K) \leq g$.  As a consequence of Theorem \ref{thm:main} we
have the following:

\begin{cor}
\label{cor:gr leq 0}

Let $M$ be a manifold, $\s \subset M$ a genus $g$ Heegaard surface
and $K \subset \s$ a knot. If $t(K) = g$, then $gr_t(K) \leq 0$.
\end{cor}

\begin{proof}

By a slight isotopy of $K \subset \s$ it is easy to see that $K$
is 1-bridge with respect to $\s$.  We can therefore apply Theorem
\ref{thm:main} to obtain the corollary.
\end{proof}

Let $T_{p,q}$ be a non-trivial torus knot.  Note that $T_{p,q}$
embeds in a genus 1 Heegaard surface of $S^3$ and $t(T_{p,q})=1$.
Hence, by Corollary \ref{cor:gr leq 0}, $gr_t(T_{p,q}) \leq 0$. On
the other hand, it is known that torus knots are meridionally
small.  Therefore by (2) of Remarks \ref{rmks:basic facts}
$gr_t(T_{p,q}) \geq 0$.  As a conclusion, we have that
$gr_t(T_{p,q}) = 0$.

Similarly, given any Heegaard surface $\s \subset M$ and a knot $K
\subset \s$ fulfilling the hypotheses of Corollary \ref{cor:gr leq
0}, by performing any Dehn surgery for which the boundary of the
meridian of the attached solid torus intersects each component of
the boundary of $\s \cap E(K)$ exactly once, we get $M'$, $\s'$
and $K'$ also fulfilling the hypotheses of Corollary \ref{cor:gr
leq 0} (see, for example, \cite{rieck}). By Hatcher \cite{hatcher}
after excluding a finite set, $K'$ is known to be meridionally
small. Therefore we obtain infinitely many knots $K'$ with
$gr_t(K') = 0$.

Suppose, particularly, that $K$ is isotopic onto some minimal
genus Heegaard surface.  Such knots are called \em good \em in
\cite{rieck} where it was shown that either $g(M) = t(K)$ or $g(M)
= t(K) + 1$. For such knots, Corollary \ref{cor:gr leq 0} implies
that if $g(M) = t(K)$, then $gr_t(K) \leq 0$.  In contrast to
this, in case when $g(M) = t(K) + 1$ (equivalently, $g(M) =
g(E(K))$) we have the following.  (Note that if $g(E(K)) = g(M)$
then $K$ is necessarily isotopic onto a minimal genus Heegaard
surface of $M$.)

\begin{thm}
\label{thm:g(M) = g(E(K))}

Let $K$ be a knot in a closed, orientable 3-manifold $M$ such that
$g(M) = g(E(K))$.  Then for all $n$,  $t(nK) = nt(K) + n-1$. In
particular, $gr_t(K) =1$, and $\lim_{m \to \infty}  \frac{t(mK) -
mt(K)}{m-1}$ exists.
\end{thm}

So we see that the growth rate is 1 if and only if $g(M) =
g(E(K))$, \ie, if and only if $K$ is a core of a handlebody in a
minimal genus Heegaard splitting of $M$.

In the rest of this section we describe a relationship between
Theorem \ref{thm:main} and Morimoto's Conjecture,  
which is concerned with the super additive phenomenon of the
tunnel number of knots. In \cite{morimoto-annalen} Morimoto
conjectured that
$t(K_1 \# K_2) \le t(K_1) + t(K_2) $  
if and only if $K_1$ or $K_2$ admits a primitive meridian (for
definitions, see Subsection \ref{subsection:primitive meridian}).
Here we note the well-known fact that if $K_1$ or $K_2$ admits a
primitive meridian, then $t(K_1 \# K_2) \le t(K_1) + t(K_2)$ (see
Proposition~\ref{pro:connecting with a PM} in section~2).

In this paper
we show that the growth rate can furnish a method for disproving
Morimoto's Conjecture. For this purpose, we first demonstrate the
following theorem with assuming Theorem~ \ref{thm:main}, and
Proposition~\ref{pro:connecting with a PM}.

\begin{thm}
\label{thm:morimoto-counterexample}

If there exists a knot $K \subset S^3$ so that neither $K$ nor $2K
(=K \# K)$ admits a primitive meridian then Morimoto's Conjecture
is false.
\end{thm}

\begin{rmk}
The following proof shows that Theorem
\ref{thm:morimoto-counterexample} still holds for knots $K \subset
M$ in an arbitrary manifold, provided that $g(E(K)) > g(M)$.
\end{rmk}

\begin{proof}
By Theorem \ref{thm:main} for some $n$ we have that

$$t(nK) < nt(K) + (n-1).$$

From now on, take $n \geq 2$ to be the minimal $n$ with that
property.

If $n=2$ then two copies of $K$ provide the desired
counterexample.  Thus we may assume that $n >2$.  Consider the knots
$K_1=2K$ and $K_2=(n-2)K$.  Note that $K_1$ does not admit a
primitive meridian by assumption and $K_2$ does not admit a
primitive meridian by Proposition \ref{pro:connecting with a PM}
and the minimality of $n$.

We have:

\begin{enumerate}
    \item $t(K_1) = 2t(K) + 1$ (by minimality of $n$),
    \item $t(K_2) = (n-2)t(K) + (n-3)$ (by minimality of $n$), and---
    \item $t(K_1 \# K_2) = t(nK) < nt(K) + (n-1) = t(K_1) + t(K_2) +1$.
\end{enumerate}

Hence
$K_1$ and $K_2$
provide a counterexample to Morimoto's Conjecture.
\end{proof}

The following corollary follows from Theorem~\ref{thm:morimoto-counterexample}, and Lemma~\ref{lem:if n exists} in section~3.

\begin{cor}
\label{cor:morimoto-counterexample}

If there exists a knot in $S^3$ with growth rate grater than $1/2$
then Morimoto's Conjecture is false.
\end{cor}

\begin{rmk}
The following proof shows that Corollary
\ref{cor:morimoto-counterexample} still holds for knots $K \subset
M$ in an arbitrary manifold, provided that $g(E(K)) > g(M)$.
\end{rmk}

\begin{proof}
Let $K \subset S^3$ be a knot with $gr_t(K) > \frac{1}{2}$. If $K$
admits a primitive meridian, then by
Proposition~\ref{pro:connecting with a PM} $t(K \# qK) \leq t(K) +
t(qK)$ for any positive integer $q$. This implies $t(nK) \le
nt(K)$ for any $n$, and consequently we have $gr_t (K) \le 0$, a
contradiction. If $K\# K$ admits a primitive meridian, then for
any $q$ we have that $t(2K \# qK) \leq t(2K) + t(qK)$, and by
Lemma \ref{lem:if n exists}, $gr_t(K) \leq \frac{1}{2}$, a
contradiction; hence neither $K$ nor $2K$ admits a primitive
meridian and  Morimoto's Conjecture is false by
Theorem~\ref{thm:morimoto-counterexample}.
\end{proof}

We conclude this section with some open questions related to the
above; in all these questions we assume that $g(M) < g(E(K))$. Let
$g$ denote the Heegaard genus of $E(K)$.

For the definitions of genus $g$, $n$-bridge presentation of a
knot, and $(g, n)$-knot in a closed 3-manifold, see
Subsection~\ref{subsection:bridge number}. It is well known
(Proposition~\ref{characterization of knots with PM from bridge
position}) that $K$ admits a primitive meridian if and only if $K$
is a $(t(K),1)$-knot, that is, $K$ has bridge index one with
respect to some genus $t(K)$ Heegaard surface of $M$. In the proof
of Theorem \ref{thm:main} we generalize this and see that:

\begin{pro}\label{prop:sect 1}
If $K$ is a $(t(K),n)$-knot then either $nK$ admits a primitive
meridian or $t(nK) < nt(K) + n-1$ (possibly both).
\end{pro}

We note that an interesting outcome of Proposition \ref{prop:sect
1} is the case when $t(nK) = nt(K) + n -1$ with $nK$ admitting a
primitive meridian (in Section \ref{sec:example} we will see such
phenomenon indeed occurs for $n=2$).  About the converse we ask:

\begin{question}
\label{que:existence of bridge presentation}

Suppose $t(nK) = nt(K) + n-1$ and that $nK$ admits a primitive
meridian. Does $K$ admit a genus $t(K)$, $n$-bridge presentation?
Specifically, suppose $t(2K) = 2t(K) + 1$ and that $2K$ admits a
primitive meridian. Does $K$ have a genus $t(K)$, 2-bridge
presentation?
\end{question}

Next we ask:

\begin{question}
\label{que:existence of high bridge number}

Does there exist a knot $K$ so that $K$ is a $(t(K),n)$ knot for
$n$ at least 3?  Specifically, does there exist a tunnel number 1
knot in $S^3$ that has bridge index greater than 2 with respect to
the genus 1 Heegaard splitting of $S^3$?
\end{question}

\begin{rmk}
If there exists a knot for which the answers to Question
\ref{que:existence of bridge presentation} (the case $n=2$) and
Question \ref{que:existence of high bridge number} are both ``yes"
then neither $K$ nor $2K$ admits a primitive meridian. Therefore,
by Theorem \ref{thm:morimoto-counterexample},
Morimoto's Conjecture is false.
\end{rmk}

\begin{question}
\label{que:spectrum}

What is the spectrum of growth rates?
\end{question}

We note that it is known that we can construct arbitrarily high
degeneration of tunnel number of knots under connected sum (see
\cite{kobayashi-degeneration}), nevertheless the following
question is still open.

\begin{question}
\label{que:negative growth rate}

Does there exist a knot with negative growth rate?  If so, can the
growth rate of knots be arbitrarily negative?
\end{question}

In light of Corollary \ref{cor:morimoto-counterexample} we ask:

\begin{question}
\label{que:knot with slope > 1/2}

Is there a knot with growth rate greater than 1/2 and less than 1?
\end{question}

We know little about the properties of the sequence $\frac{t(mK) -
mt(K)}{m-1}$.  In particular:

\begin{question}
\label{que:properties of sequence}

Is the sequence $\frac{t(mK) - mt(K)}{m-1}$ eventually monotonous?

Does $\lim_{m \to \infty}  \frac{t(mK) - mt(K)}{m-1}$ exist?
\end{question}

\begin{question}
Given a knot $K$, is $t(nK) \leq t((n+1)K)$?
\end{question}

We remark that while a positive answer would give bounds on the
behavior of the elements $\frac{t(nK) - nt(K)}{n-1}$, it will not
suffice to show that the limit exists.

{\bf Acknowledgements:} We thank Chaim Goodman--Strauss and Mark
Johnson for helpful conversations.

\section{Preliminaries}
\label{sec:prelims}

Throughout this paper we work in the smooth category. We always
assume our manifold to be compact and orientable.  For standard
notion in 3-manifold topology we refer the reader to \cite{hempel}
or \cite{jaco}.

\subsection{Amalgamation of Heegaard splittings}
\label{subsection:amalgamation}

A 3-manifold $C$ is a \em compression body \em if there is a
closed connected surface $F$ such that $C$ is obtained from
$F\times [0,1]$ by attaching 2-handles along mutually disjoint
\scc s in $F \times \{1\}$ and capping off the resulting 2-sphere
boundary components which are disjoint from $F \times \{0\}$ by
3-handles. The boundary component of $C$ corresponding to $F
\times \{0\}$ is denoted $\del_+ C$.  Then $\del_- C = \del C
\setminus \del_+ C$. A compression body $C$ is called a \em
handlebody \em if $\del_- C = \emptyset$.  A compressing disk $D
\subset C$ for $\del_+ C$ is called a \em meridian disk \em of the
compression body $C$.

By extending the cores of the 2-handles in the definition of the
compression body $C$ vertically to $F \times [0,1]$ we obtain a
collection of mutually disjoint meridian disks of $C$, say
$\mathcal{D}$, such that the manifold obtained from $C$ by cutting
along $\mathcal{D}$ is homeomorphic to a union of $\del_- C \times
[0,1]$ and a (possibly empty) collection of 3-balls.  This gives a
dual description of the compression body, that is, a connected
3-manifold $C$ is a compression body if there exists a (not
necessarily connected and possibly empty) closed surface
$\mathcal{F}$, without 2-sphere components, and a
(possibly empty) collection of 3-balls $\mathcal{B}$ such that $C$
is obtained from $\mathcal{F} \times [0,1] \cup \mathcal{B}$ by
attaching 1-handles to $\mathcal{F} \times \{0\} \cup \del
\mathcal{B}$.  We note that $\del_- C$ is the surface
corresponding to $\mathcal{F} \times \{1\}$.

Let $N$ be a compact 3-manifold and $F_1$, $F_2$ a partition of
the components of $\del N$.  We say that a decomposition $C_1
\cup_{\s} C_2$ (or $C_1 \cup C_2$) is a Heegaard splitting of
$(N;F_1,F_2)$ (or of $N$) if it satisfies the following
conditions:

\begin{enumerate}
    \item $C_i$ ($i=1,2$) is a compression body in $N$ such that
    $\del_- C_i = F_i$,
    \item $C_1 \cup C_2 = N$, and---
    \item $C_1 \cap C_2 = \del_+ C_1 = \del_+ C_2 = \s$.
\end{enumerate}

The surface $\s$ is called a \em Heegaard surface \em of
$(N;F_1,F_2)$ (or $N$).  The genus of $\s$ is called the genus of
the splitting and the genus of the minimal genus Heegaard
splitting for $N$ is called the \em Heegaard genus of $N$, \em
denoted $g(N)$.

Let $\s_1 \subset M_1$ and $\s_2 \subset M_2$ be two Heegaard
surfaces for 3-manifolds $M_1$ and $M_2$.  Suppose that a
component $T_1$ of $\del M_1$ and a component $T_2$ of $\del M_2$
are homeomorphic. Let $M$ be a manifold obtained from $M_1$ and
$M_2$ by identifying $T_1$ and $T_2$ by a homeomorphism.  Then we
can obtain a Heegaard surface for $M$, say $\s$, from $\s_1$ and
$\s_2$ by collapsing the product region adjacent to $T_1$ and
$T_2$ as in Figure~1. We call $\Sigma$ the \em amalgamation \em of
$\s_1$ and $\s_2$ along $T$, where $T$ is the image of $T_1 = T_2$
in $M$. For details, see \cite{schultens-FXS1}, where it was shown
that:

\begin{equation}
\label{equ:genus of amalgamation}
g(\s) = g(\s_1) + g(\s_2) - g(T).
\end{equation}

\begin{figure}[ht]
\begin{center}
\includegraphics[width=10cm, clip]{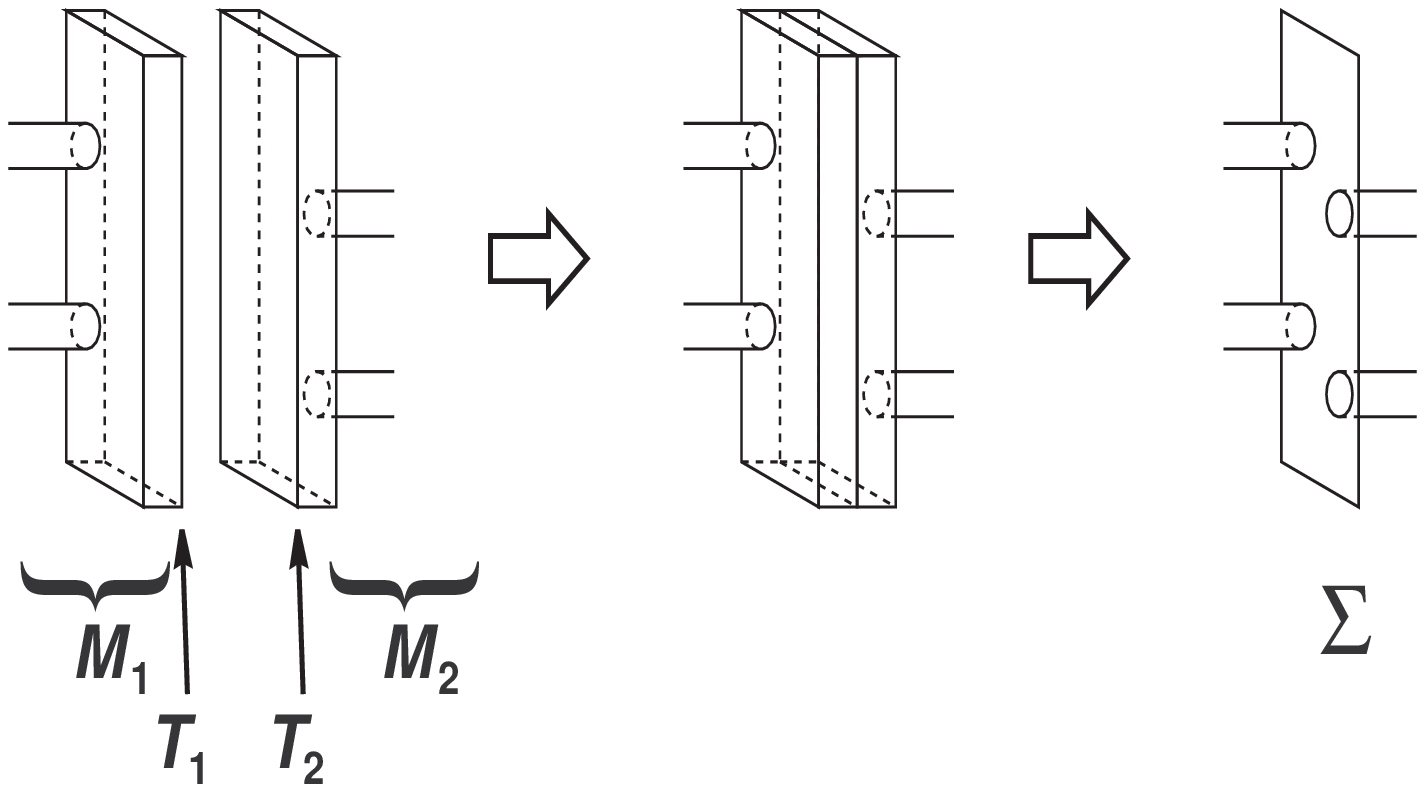}
\end{center}

\begin{center}
Figure 1
\end{center}

\end{figure}

\subsection{Tunnel number}
\label{subsection:tunnel}

Let $M$ be a compact orientable 3-manifold.  By a knot $K$ we mean
a smooth embedding of $S^1$ into $M$.  For a knot $K \subset M$,
let $E(K)$ denote the exterior of $K$, \ie, $E(K)= M \setminus
int(N(K))$.  It is an easy consequence of the existence of smooth
structures on $E(K)$ that there exists a collection of properly
embedded arcs $\tau$ such that $\text{cl}(E(K)\setminus N(\tau))$
is a handlebody. Such a collection $\tau$ is called a \em tunnel
system. \em The \em tunnel number \em of $K$, denoted $t(K)$, is
the minimal number of arcs required for tunnel systems.  Then
$t(K)$ is related to the Heegaard genus of $E(K)$ by the equation
$t(K) = g(E(K))-1$.  Given two knots $K_1 \subset M_1$ and $K_2
\subset M_2$ the connected sum of $K_1$ and $K_2$, denoted $K_1\#
K_2$, is a knot in $M_1 \# M_2$ which is obtained by removing from
$M_i$ a small ball so that $K_i$ intersects this ball 
in a single unknotted arc ($i=1,2$) and then identifying the
boundaries of the punctured manifolds by a homeomorphism under
which the intersection of the knots with the boundaries match up.
By taking a union of a tunnel system for $K_1$, a tunnel system
for $K_2$ and one extra tunnel on the decomposing sphere we obtain
a tunnel system for $K_1 \# K_2$.
%
This gives the Inequality~(\ref{eq:trivial}) of Section~\ref{sec:intro}.

When equality in Inequality~(\ref{eq:trivial}) holds
we say that the tunnel number of the knots
are \em super additive. \em

\subsection{Primitive meridian}
\label{subsection:primitive meridian}

Let $M$, $K$ be as above.  We say that $K$ admits a \em primitive
meridian \em if there is a minimal genus Heegaard splitting of
$E(K)$, say $C \cup_{\s} V$ (here $C$ is a compression body and
$V$ is a handle body, so that $\del E(K) = \del_- C$), such that
there exists a properly embedded essential annulus $A \subset C$
and a meridian disk $D \subset V$ with the following properties:

\begin{enumerate}
    \item $A$ is vertical in $C$ (\ie, one boundary component of
    $A$ is on $\del_+ C$ and the other on $\del_- C$) and $A \cap
    \del_- E(K)$ is a meridian curve of $E(K)$.
    \item $A \cap \s$ ($=A \cap \del_+ C$) and $D \cap
    \s$ ($=\del D$) intersect transversely (in $\s$) in one point.
\end{enumerate}

An important feature of knots admitting a primitive meridian is
that they are never super additive (recall Subsection
\ref{subsection:tunnel}) when connected sum to other knots. In
fact, Propositions 1.3 and 2.1 of \cite{morimoto-annalen} imply
the following:

\begin{pro}
\label{pro:connecting with a PM}

If $K$ admits a primitive meridian then for any knot $K'$ we have:

$$t(K \# K') \leq t(K) + t(K').$$
\end{pro}

In \cite{morimoto-annalen} Morimoto conjectures that the converse
of this is true:

\begin{cnj}[Morimoto's Conjecture]
\label{cnj:morimoto}

If $t(K \# K') \leq t(K) + t(K')$, then $K$ or $K'$ admits a
primitive meridian.
\end{cnj}

\begin{rmk}
This conjecture is stated in a different appearance in
\cite{morimoto-annalen}.  By Proposition 2.1 of that paper (see
Proposition~\ref{characterization of knots with PM from bridge
position} below), we can show that the two are equivalent.
\end{rmk}

In \cite{kobayashi-rieck} the authors give a necessary and
sufficient condition for a knot $K$ to admit a primitive meridian.
We introduce the result here.

Notation: For a knot $K$, we denote by $\widehat{K}$ the link
obtained from $K$ by adding a single \scc\ parallel to the
meridian of $K$; in other words, $\widehat{K} = K \# \mbox{ (Hopf
link in }S^3)$.  For $n \geq 0$ we denote by $K^{(n)}$ the link
obtained from $K$ by adding $n$ \scc s simultaneously parallel to
$n$ disjoint copies of the meridian in $E(K)$ (here we understand
$K^{(0)} = K$ and $K^{(1)} = \widehat{K}$).  The exterior of
$\widehat{K}$ is denoted by $\widehat{E(K)}$ and the exterior of
$K^{(n)}$ is denoted by $E(K)^{(n)}$.  An essential annulus $A$
properly embedded in $\widehat{E(K)}$ (resp. $E(K)^{(n)}$) is
called a \em Hopf spanning annulus \em if one boundary component
of $A$ is a meridian of $K$ and the other a longitudinal curve of
$\del \widehat{E(K)} \setminus \del E(K)$ (resp. some component of
$\del E(K)^{(n)} \setminus \del E(K)$).

Then the following holds:

\begin{pro}
\label{pro:characterization of knots with PM}

A knot $K \subset M$ admits a primitive meridian if and only if
there exists a Heegaard splitting $U \cup_{\s} V$ of
$\widehat{E(K)}$ which satisfies the following conditions:

\begin{enumerate}
    \item The genus of $\s$ equals $g(E(K))$,
    \item There exists a Hopf spanning annulus $A \subset
    \widehat{E(K)}$ such that $A$ intersects $\s$ transversely in a
    single \scc\ that is essential in $A$.
\end{enumerate}
\end{pro}

The proof of this proposition can be found in
\cite{kobayashi-rieck}. In this paper we only use the ``if" part
of this proposition and the idea of its proof is as follows: since
$E(K)$ is obtained form $\widehat{E(K)}$ by Dehn filling, any
Heegaard surface for $\widehat{E(K)}$ is also a Heegaard surface
for $E(K)$.  The Heegaard surface for $\widehat{E(K)}$ that is
described in the proposition naturally admits a primitive
meridian.  See Figure~2.

\begin{figure}[ht]
\begin{center}
\includegraphics[width=10cm, clip]{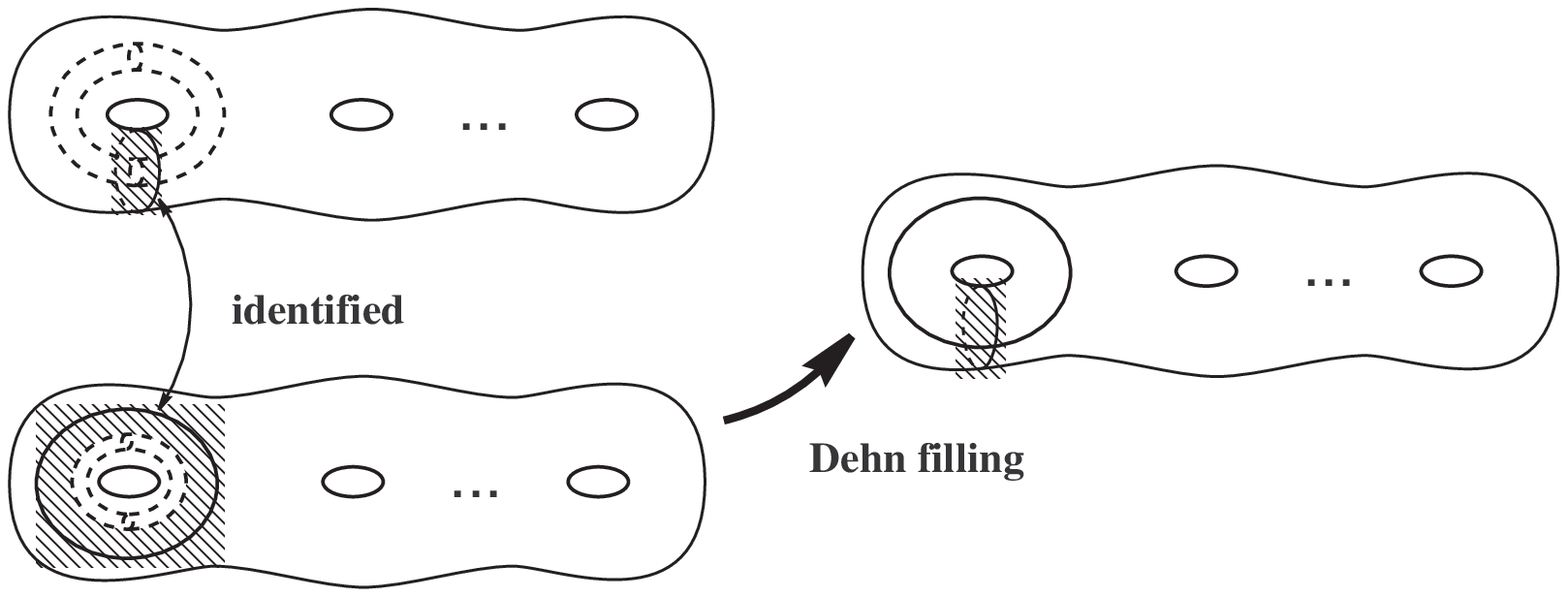}
\end{center}

\begin{center}
Figure 2
\end{center}

\end{figure}

\subsection{Generalized bridge number}
\label{subsection:bridge number}

Let $U \cup_\s V$ be a Heegaard splitting of some 3-manifold $M$,
and $K \subset M$ a knot.  Following Doll \cite{doll} we say that
$K$ is in \em genus $g$, $n$-bridge presentation \em (with respect
to $\s$) if $K\cap U$ ($K \cap V$ resp.) consists of $n$ arcs, and
these arcs are simultaneously parallel into $\del U$ ($\del V$
resp.) (The genus is omitted when clear from context). We say that
$K$ is a \em $(g,n)$-knot\em, or a genus $g$, $n$-bridge knot, if
$K$ admits a genus $g$, $n$ bridge presentation but does not admit
a genus $g$, $n-1$ bridge presentation for any Heegaard surface of
genus $g$.

Then we have the following characterization of knots with
primitive meridian (see Proposition 2.1 of
\cite{morimoto-annalen}):

\begin{pro}
\label{characterization of knots with PM from bridge position}

Let $K$, $M$ be as above.  Then $K$ admits a primitive meridian if
and only if $K$ is a $(t(K),1)$-knot (or, $(g(E(K))-1,1)$-knot).
\end{pro}

\section{The proof}
\label{sec:proof}

In this section we prove Theorems \ref{thm:main}, and
\ref{thm:g(M) = g(E(K))}. For the proof of Theorem \ref{thm:main},
we first prepare two lemmas:

\begin{lem}
\label{lem:when t goes down in nK} Let $K$ be a knot in a closed
orientable 3-manifold $M$.
Suppose there exists a positive number $n$ so that $t(nK) < nt(K)
+ n - 1$.
Then $gr_t(K) \leq \frac{n-1}{n}$.
\end{lem}

\begin{proof}

Fix a non-negative integer $m$.  Write $m$ as $pn + l$, with $0
\leq l < n$.  We have inequalities below. (Notes: (1) In the
following we take $t(0K)$ to be $0$. (2) In one line of the following lines we need
to treat the cases $l=0$, $l>0$ separately. We can deal with both
cases by introducing a variable $q$ such that $q=p$ if $l>0$, and
$q=p-1$ if $l=0$. (3) Some lines contain explanations in
[brackets]; there by Equation~(\ref{eq:trivial}) we mean
Equation~(\ref{eq:trivial}) of Section~\ref{sec:intro}.)

\begin{eqnarray*}
  t(mK) &=&  t((np+l)K)\\
  &\leq& t(nK) + t([(p-1)n + l]K) +1
                    \hspace{.6in}\mbox{[By Equation (\ref{eq:trivial})]}\\
  &\leq& 2t(nK) + t([(p-2)n + l]K) +2
                    \hspace{.5in}\mbox{[By Equation (\ref{eq:trivial})]}\\
  &\cdots \\
  &\leq& pt(nK) + t(lK) + q \\
  &\leq& p(t(nK)) + l t(K) + (l-1) +p \\
  &\leq& p(nt(K) + n-2) + l t(K) + (l-1) +p \\
  &     &  \hspace{.7in}[\mbox{since }t(nK) < nt(K) + n - 1\mbox{ by assumption}]\\
  &\leq& (np+l) t(K) + (np+l-1) - p \\
  &=& mt(K) + (m-1) - p.
\end{eqnarray*}

We conclude that the growth rate fulfills the inequality
$\frac{t(mK) - mt(K)}{m-1} \leq \frac{m-1 -p}{m-1}$. In the limit
$m \approx np$ and we get the bound $\frac{(n-1)p}{np} =
\frac{(n-1)}{n}$ as desired.  This completes the proof of Lemma
\ref{lem:when t goes down in nK}.
\end{proof}

\begin{lem}
\label{lem:if n exists}
Let $K$ be a knot in a closed orientable 3-manifold $M$.
Suppose that there exists a positive integer $n$ so that for any
positive integer $q$ we have that $t((n+q)K) \leq t(nK) + t(qK)$.
Then $gr_t(K) \leq \frac{n-1}{n}$.
\end{lem}

\begin{proof}

Again let $m$ be a positive integer and write $m$ as $pn + l$,
with $0 \leq l < n$. By repeatedly applying the assumption of
Lemma \ref{lem:if n exists}, and then applying
Equation~(\ref{eq:trivials}) in section~1 twice we have: (Note:
When we apply Equation~(\ref{eq:trivials}) to $t(lK)$ there are
two possible results, depending on whether $l=0$ or not; these are
treated by using $\max\{l-1,0\}$.  Eventually we arrive at a bound
that is valid in both cases.) 

\begin{eqnarray*}
  t(mK) &=&  t((np+l)K)\\
  &\leq& t(nK) + t([(p-1)n + l]K) \hspace{.7in} [\mbox{by
                                                  assumption}]\\
  &\leq& 2t(nK) + t([(p-2)n + l]K)  \hspace{.5in} [\mbox{by
                                                  assumption}]\\
  &\cdots \\
  &\leq& pt(nK) + t(lK)  \hspace{1.3in} [\mbox{by
                                                  assumption}]\\
  &\leq& p[nt(K) + (n-1)] + t(lK) \hspace{.6in}
                                [\mbox{by Equation~(\ref{eq:trivials})}]\\
  &\leq& p[nt(K) + (n-1)] + l t(K) + \max{\{l-1,0\}}\\
  &\leq& (np+l) t(K) + (np+ \max{\{l-1,0\}}) - p\\
  &\leq& mt(K) + m - p.
\end{eqnarray*}

Note that if $l \neq 0$ we actually get the bound $mt(K) + (m-1) -
p$, which is exactly the bound we obtained in the proof of Lemma
\ref{lem:when t goes down in nK}.  Since the constant +1 vanishes
in the limit, we complete the proof of Lemma \ref{lem:if n exists}
exactly as above.
\end{proof}

\begin{proof}[Proof of Theorem \ref{thm:main} and Proposition \ref{prop:sect 1}]

Let $K \subset M$, $n$ be as in Theorem \ref{thm:main}.  Suppose
that $n=1$.
By Proposition~\ref{characterization of knots with PM from bridge position},
$K$ admits a primitive meridian.
Hence, by applying Proposition \ref{pro:connecting with a PM}
repeatedly, we see that for any positive integer $m$ we have that
$t(mK) \leq mt(K)$.  Hence $gr_t(K) \leq 0 = \frac{1-1}{1}$; this
gives Theorem \ref{thm:main}.  Hence for the remainder of the
proof we assume that $n \geq 2$.

If $t(nK) < nt(K) +(n-1)$ then by Lemma \ref{lem:when t goes down
in nK} we have the conclusion of Theorem \ref{thm:main}.  Hence
for the remainder of the proof we may assume that $t(nK) = nt(K)
+(n-1)$, and will show that $nK$ admits a primitive meridian.
(Note that this proves Proposition~\ref{prop:sect 1}.) Let $U
\cup_\s V$ be a genus $t(K)$ Heegaard splitting of $M$ with
respect to which $K$ admits an $n$ bridge presentation. Isotope
$K$ to such position.  Hence $K \cap U$ ($K \cap V$ resp.)
consists of $n$ boundary parallel arcs.  Then $E(K) \cap U$ ($E(K)
\cap V$ resp.) is a genus $t(K) + n$ handlebody, say $U'$ ($V'$
resp.) such that $U' \cap \del E(K)$ ($V' \cap \del E(K)$ resp.)
is a collection of $n$ disjoint simultaneously boundary
compressible annuli in $\del U'$ ($\del V'$ resp.) as in Figure~3;
these annuli are meridional in $\del E(K)$.

\begin{figure}[ht]
\begin{center}
\includegraphics[width=10cm, clip]{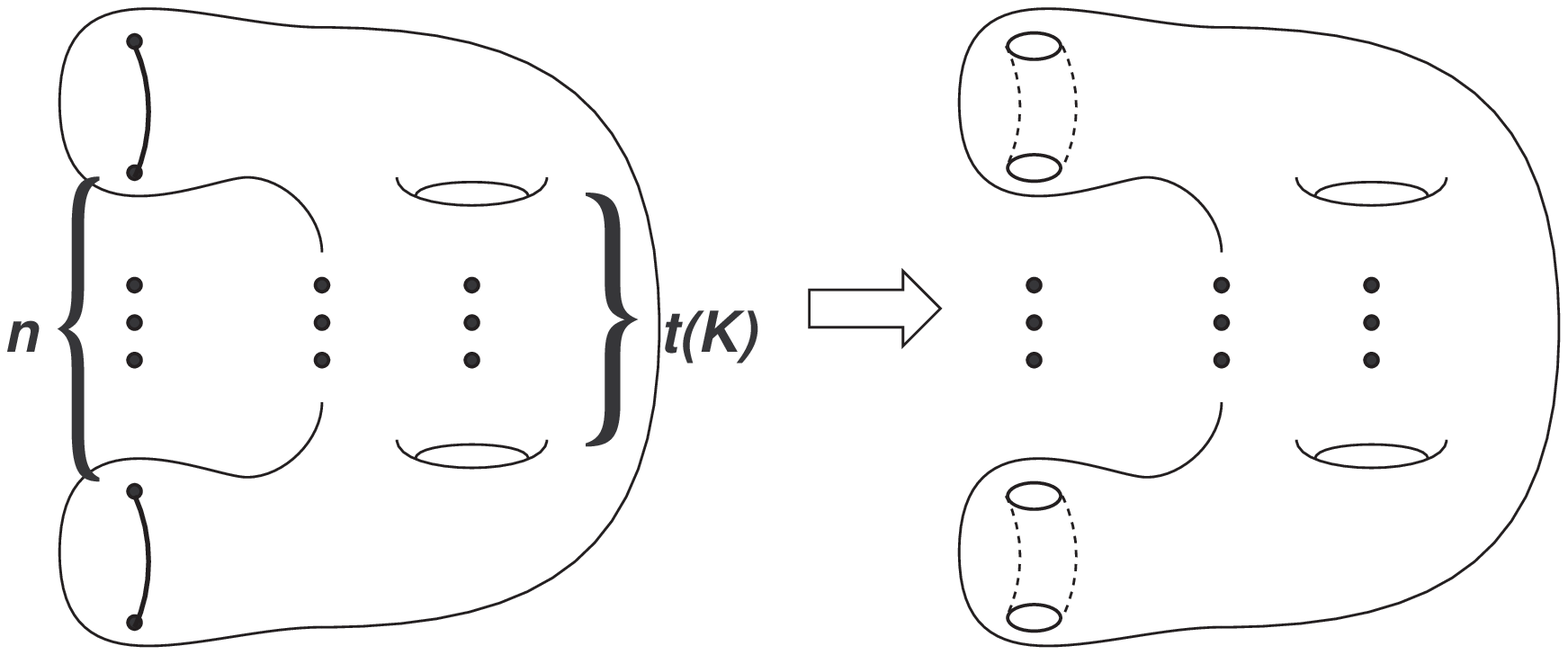}
\end{center}

\begin{center}
Figure 3
\end{center}

\end{figure}

Recall from Subsection \ref{subsection:primitive meridian} the
notation $K^{(n)}$.

\medskip
\noindent {\bf Claim.} $E(K)^{(n)}$ admits a genus $t(K) + n$
Heegaard splitting, say $\widetilde{U} \cup_{\widetilde{\s}}
\widetilde{V}$, such that $\widetilde{\s}$ intersects a Hopf
spanning annulus in a single \scc\ which is essential in the
annulus.

\medskip
\begin{proof}
Let $A^{(n)}$ be an annulus with $n$ holes.  Then let $l^+$, $l^-$
be the boundary components of the annulus, and let
$l_1,\cdots,l_n$ be the boundary components corresponding to the punctures.

\begin{figure}[ht]
\begin{center}
\includegraphics[width=6cm, clip]{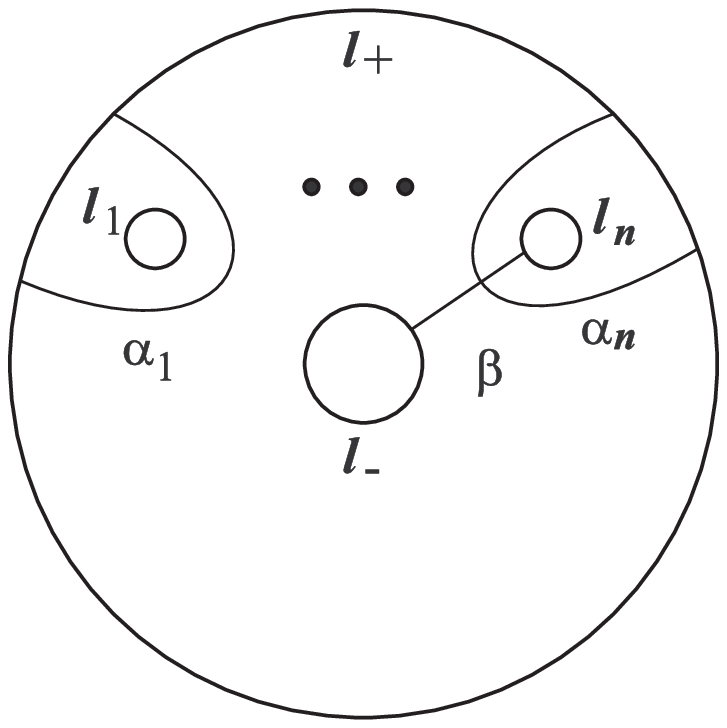}
\end{center}

\begin{center}
Figure 4
\end{center}

\end{figure}

We note that $E(K)^{(n)}$ can be represented as a manifold
obtained from $E(K)$ by attaching $A^{(n)} \times S^1$ by a
homeomorphism $h: \del E(K) \to l_+ \times S^1$ that maps a
meridian curve to $\{pt\} \times S^1$. Let
$\alpha_1,\cdots,\alpha_n$ be arcs properly embedded in $A^{(n)}$
as in Figure~4.  Note that the annuli $\cup_{i=1}^n (\alpha_i
\times S^1)$ cut $A^{(n)}$ into $n+1$ pieces, say $T_0,
T_1,\cdots,T_n$ each of which is homeomorphic to (annulus)$\times
S^1$ (or (torus)$\times [0, 1]$), where $l_- \times S^1 \subset
T_0$ and for each $i$, $l_i \times S^1 \subset T_i$.   Recall that
$\partial E(K) \cap U'$ ($\partial E(K) \cap V'$ resp.) consists
of $n$ annuli. Since $h$ maps a meridian curve to $\{pt\} \times
S^1 \subset l_+ \times S^1$, we may suppose, by deforming $h$ by
an isotopy if necessary, that $h(E(K) \cap U') = \cup_{i=1}^n (T_i
\cap (l_+ \times S^1))$ (hence, $h(\del E(K) \cap V')= T_0 \cap
(l_+ \times S^1)$).  Then let $\widetilde{U} = U' \cup_h
(\cup_{i=1}^n T_i)$ and $\widetilde{V} = V' \cup_h T_0$.

\begin{figure}[ht]
\begin{center}
\includegraphics[width=10cm, clip]{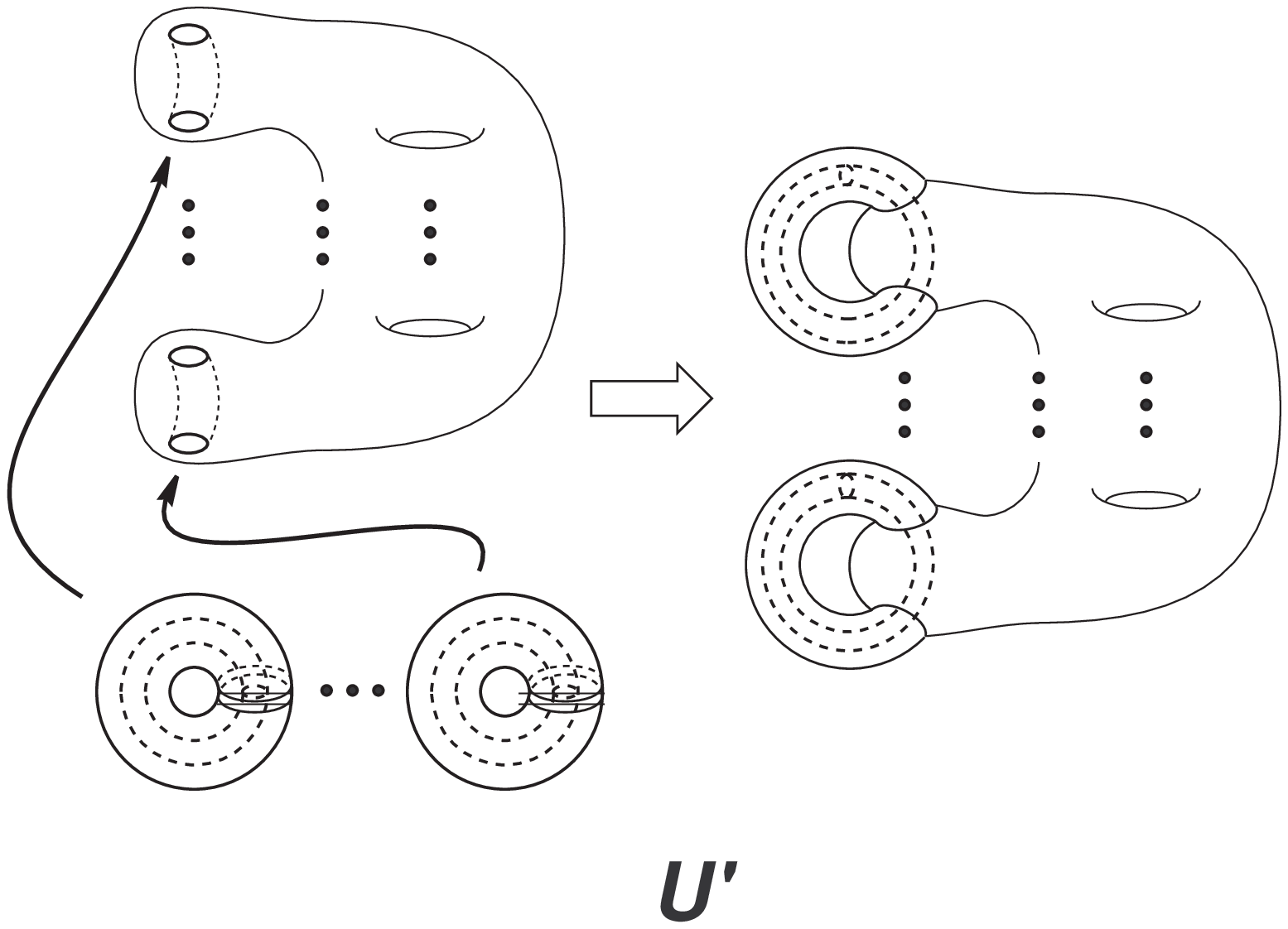}
\end{center}

\begin{center}
Figure 5
\end{center}

\end{figure}

\begin{figure}[ht]
\begin{center}
\includegraphics[width=10cm, clip]{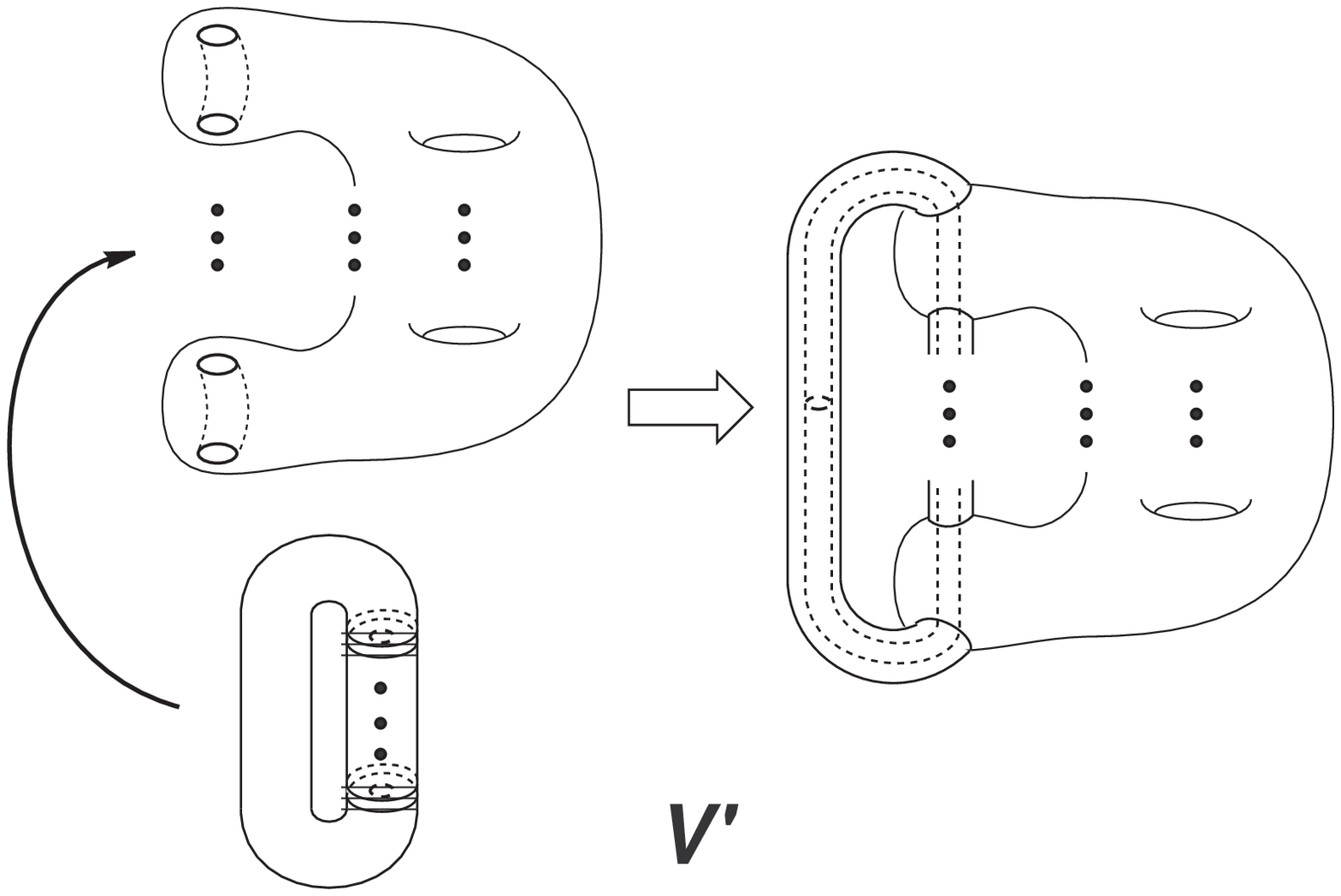}
\end{center}

\begin{center}
Figure 6
\end{center}

\end{figure}

By Figures~5 and 6 we have that $\widetilde{U}$ and
$\widetilde{V}$ are compression bodies of genus $t(K) + n$, and
$\widetilde{U} \cap \widetilde{V} = \del_+ \widetilde{U} = \del_+
\widetilde{V} (=\widetilde{\s})$.  Hence $\widetilde{U}
\cup_{\widetilde{\s}} \widetilde{V}$ is a genus $t(K) + n$
Heegaard splitting of $E(K)^{(n)}$.  Let $\beta$ be the arc
properly embedded in $A^{(n)}$ as in Figure~4.  Note that $\beta
\times S^1$ gives a Hopf spanning annulus in $E(K)^{(n)}$ and it
is directly observed that it intersects the Heegaard surface
$\widetilde{\s}$ in a single \scc\ corresponding to $(\alpha_n
\cap \beta) \times S^1$. This gives the conclusion of the claim.
\end{proof}

Note that $\widehat{E(nK)}$ is obtained from $E(K)^{(n)}$ and
$n-1$ copies of $E(K)$ by identifying $l_1 \times
S^1,\cdots,l_{n-1} \times S^1$ and the boundary components of the
copies of $E(K)$ by a homeomorphism which takes each $\{ pt \}
\times S^1 \subset l_i \times S^1$ to a meridian curve. Let $Q$ be
the Heegaard surface of $\widehat{E(nK)}$ obtained from
$\widetilde{\s}$ and a minimal genus Heegaard surface for each of
the copies of $E(K)$ by amalgamations. By Equation~(\ref{equ:genus
of amalgamation}) of Subsection \ref{subsection:amalgamation}, the
genus of $Q$ fulfills:

\begin{eqnarray*}
  g(Q) &=& (t(K) + n) + (n-1) g(E(K)) - (n-1) \\
   &=& t(K) + n + (n-1) (g(E(K))-1) \\
   &=& t(K) + n+ (n-1)t(K)\\
   &=& n(t(K)+1) \\
   &=& n g(E(K)).
\end{eqnarray*}

Recall that we assumed that $t(nK) = nt(K) + (n-1)$ (equivalently
$g(E(nK)) = ng(E(K))$) at the beginning of this proof.  Hence we
have $g(Q) = g(E(nK))$.  On the other hand $Q$ is a Heegaard
surface for $\widehat{E(nK)}$; hence $g(\widehat{E(nK)}) \leq
g(E(nK))$. Since every Heegaard surface of $\widehat{E(nK)}$ is
naturally a Heegaard surface of $E(nK)$, we see that
$g(\widehat{E(nK)}) \geq g(E(nK))$.  Hence we have
$g(\widehat{E(nK)}) = g(E(nK))$. The amalgamation above does not
effect the annulus $\beta \times S^1$ and its intersection with
$\widetilde{\s}$.  Hence by Proposition \ref{pro:characterization
of knots with PM} we see that $nK$ admits a primitive meridian.
Then by Proposition \ref{pro:connecting with a PM} we see that
$t((n+q)K) \le t(nK) + t(qK)$ for any $q$. Then by Lemma
\ref{lem:if n exists} we see that $gr_t(K) \leq \frac{n-1}{n}$.

This completes the proof of Theorem \ref{thm:main}.
\end{proof}

\begin{rmk}

The arguments in the proof of Theorem~\ref{thm:main} work for a
bridge number $n$ with respect to a Heegaard surface of genus
$t(K)-k$ for $k>0$. However, it is elementary to show (see, for
example, Figure~4 of \cite{kobayashi-2-bridge}) that we can obtain
a $(t(K), n-k)$ bridge presentation from the $n$-bridge
presentation. Hence we can show that $gr_t(K) \leq
\frac{n-k-1}{n-k}$ by Theorem~\ref{thm:main}. Note that
$\frac{n-k-1}{n-k} < \frac{n-1}{n}$. This shows that using our
techniques, the best estimate of the growth rate is obtained by
using the bridge index with respect to genus $t(K)$ Heegaard
surfaces.
\end{rmk}

For the proof of Theorem~\ref{thm:g(M) = g(E(K))}, we prepare
three lemmas.

Let $\ell$ be a knot in a handlebody $V$.
We say that $\ell$ is {\it primitive} in $V$ if
there is a meridian disk $D$ of $V$ such that $D$ cuts off a solid torus,
such that $\ell$ is a core curve of the solid torus.

\begin{lem}
\label{lem:lem 1 for g(M) = g(E(K))} Let $K$ be a knot in a
closed, orientable 3-manifold $M$. Then the equality $g(M) =
g(E(K))$ holds if and only if there exists a minimal genus
Heegaard splitting $V \cup_{\Sigma} W$ of $M$ such that $K \subset
V$, and $K$ is primitive in $V$.
\end{lem}

The proof is easy, and we omit it.

\begin{lem}
\label{lem:lem 2 for g(M) = g(E(K))} Let $K$ be a knot in a
closed, orientable 3-manifold $M$. Then the equality $g(M) =
g(E(K))$ holds if and only if there exists a minimal genus
Heegaard splitting $V \cup_{\Sigma} W$ of $M$ such that $K \subset
\Sigma$, and there is a meridian disk $D$ of $V$ such that
$\partial D$ and $K$ intersects transversely (in $\Sigma$) in one point.
\end{lem}

\begin{proof}
Only if part:
Suppose $g(M) = g(E(K))$.
By Lemma~\ref{lem:lem 1 for g(M) = g(E(K))},
there exists a minimal genus Heegaard splitting $V \cup_{\Sigma} W$
of $M$ such that $K \subset V$, and
$K$ is primitive in $V$,
i.e.,
there is a meridian disk $D'$ of $V$ such that $D'$ cuts off a solid
torus $T$ such that $K$ is a core curve of $T$.
Then we can isotope $K$ so that $K \subset \partial T$,
and $K \cap D' = \emptyset$ (hence $K \subset \Sigma$).
Note that $K$ is a longitude of $T$
and this shows that we take a meridian disk $D$ of $V$ such that
$\partial D$ and $K$ intersects transversely in one point
and that $D \cap D' = \emptyset$
(hence $D$ is a meridian disk of $V$).

\medskip
If part:
Suppose that $K \subset \Sigma$ and there is a meridian disk $D$ of $V$
such that $\partial D$ and $K$ intersects transversely in one point.
Let $N(K \cup D, V)$ be a regular neighborhood of $K \cup D$ in $V$,
and $D'$ the frontier of $N(K \cup D, V)$ in $V$.
It is directly observed that $N(K \cup D, V)$ is a solid torus
such that $K$ is a longitude of $N(K \cup D, V)$, and that
$D'$ is a meridian disk of $V$ which cuts off $N(K \cup D, V)$.
Note that $K$ can be isotoped to a core curve of $N(K \cup D, V)$
by an isotopy not affecting $D'$.
Hence by Lemma~\ref{lem:lem 1 for g(M) = g(E(K))},
we see that $g(M) = g(E(K))$.

This completes the proof of Lemma~\ref{lem:lem 2 for g(M) = g(E(K))}
\end{proof}

\begin{lem}
\label{lem:lem 3 for g(M) = g(E(K))} Let $K_i$ $(i=1,2)$ be knots
in a closed, orientable 3-manifolds $M_i$. Suppose that $g(M_1) =
g(E(K_1))$, and $g(M_2) = g(E(K_2))$ hold. Then we have $g(M_1 \#
M_2) = g(E(K_1 \# K_2))$.
\end{lem}

\begin{proof}
By Lemma~\ref{lem:lem 2 for g(M) = g(E(K))} there exist minimal
genus Heegaard splittings $V_i \cup_{\Sigma_i} W_i$ of $M_i$
$(i=1,2)$ such that $K_i \subset \Sigma_i$, and there is a
meridian disk $D_i$ of $V_i$ such that $\partial D_i$ and $K_i$
intersects transversely in one point. Then take a 3-ball $B_i$ in
$M_i$ such that $\Sigma_i \cap B_i$ is a disk properly embedded in
$B_i$, $K_i \cap B_i$ is an arc properly embedded in $\Sigma_i
\cap B_i$, and $B_i \cap D_i = \emptyset$. Then we identify the
boundaries of $\text{cl}(M_1 \setminus B_1)$, and $\text{cl}(M_2
\setminus B_2)$ by a homeomorphism which identifies $\partial
(\Sigma_1 \setminus B_1)$ with $\partial (\Sigma_2 \setminus B_2)$
and $\partial (K_1 \setminus B_1)$ with $\partial (K_2 \setminus
B_2)$ to obtain a connected sum $K_1 \# K_2$. We note that the
image of $\text{cl}(\Sigma_1 \setminus B_1) \cup
\text{cl}(\Sigma_2 \setminus B_2)$ is a minimal genus Heegaard
surface of $M_1 \# M_2$ (Proposition II.10 of \cite{jaco}), and
$K_1 \# K_2$ is contained in the Heegaard surface. Note also that
$D_1$ survives in the Heegaard splitting as a meridian disk, and
$\partial D_1$ intersects $K_1 \# K_2$ transversely in one point.
Hence by Lemma~\ref{lem:lem 2 for g(M) = g(E(K))}, we see that
$g(M_1 \# M_2) = g(E(K_1 \# K_2))$.

This completes the proof of Lemma~\ref{lem:lem 3 for g(M) = g(E(K))}.
\end{proof}

\begin{proof}[Proof of Theorem \ref{thm:g(M) = g(E(K))}]
Let $K \subset M$, $n$ be as in Theorem \ref{thm:g(M) = g(E(K))}.
Since $g(M) = g(E(K))$, we have $g(E(nK)) = g( \#_{i=1}^n M)$ by
Lemma~\ref{lem:lem 3 for g(M) = g(E(K))}. Note that $g( \#_{i=1}^n
M) = n g(M)$ by Haken (Proposition II.10 of \cite{jaco}). Hence we
have:

\begin{eqnarray*}
  t(nK)&=& g(E(nK)) -1 \\
   &=& n g(M) -1 \\
   &=& n (t(K) +1) -1 \\
   &=& n t(K) + (n-1).
\end{eqnarray*}

This completes the proof of Theorem~\ref{thm:g(M) = g(E(K))}.
\end{proof}

\section{Example}
\label{sec:example}

Let $K_m = K(7,17;10m-4)$ be the knot introduced by
Morimoto--Sakuma--Yokota in \cite{skuma-morimoto-yokata}. That is,
$K$ is obtained from the torus knot $K(7,17)$ by adding twists
along an unknotted curve $\gamma$ as in Figure~7.

\begin{figure}[ht]
\begin{center}
\includegraphics[width=10cm, clip]{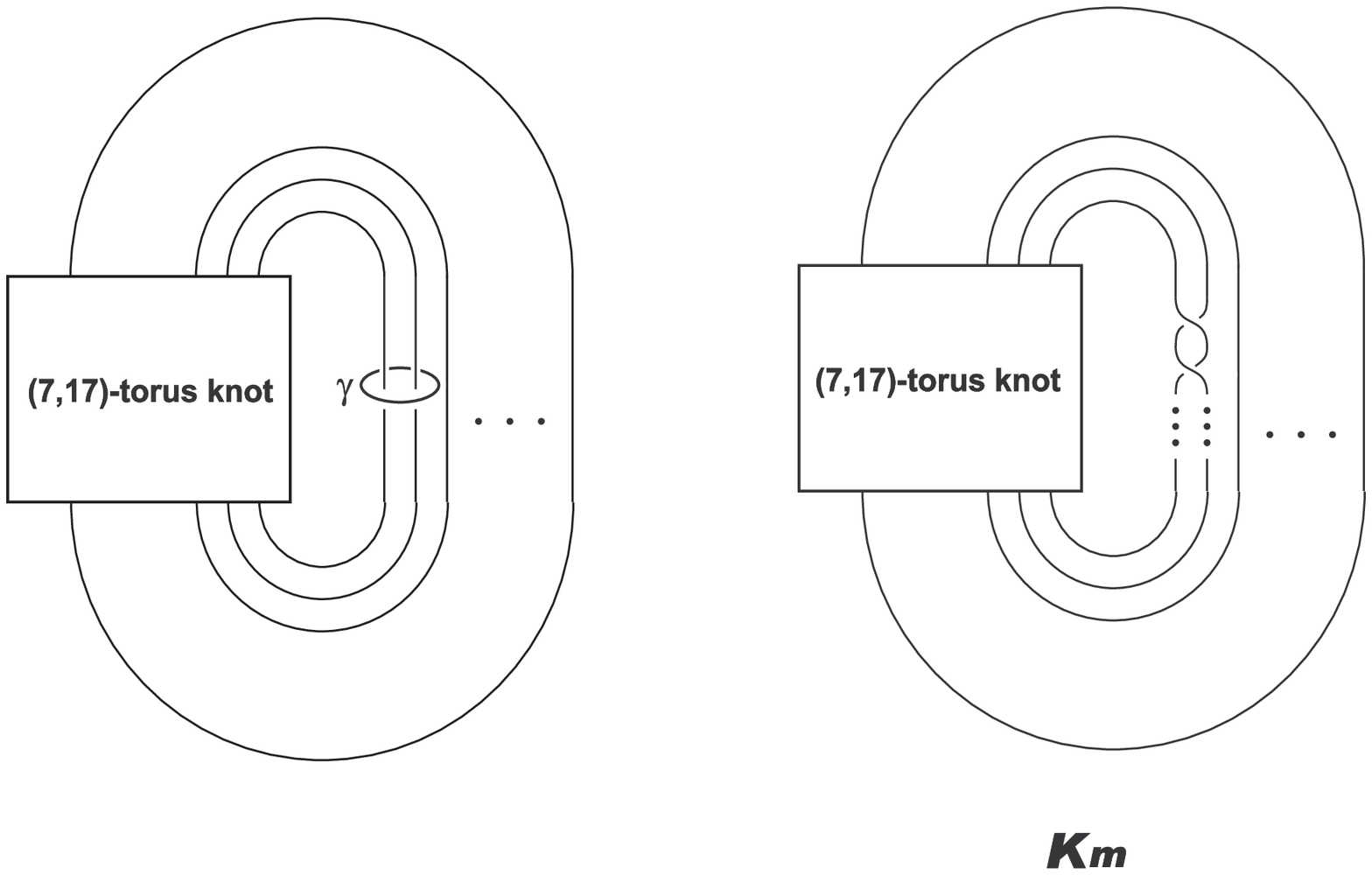}
\end{center}

\begin{center}
Figure 7
\end{center}

\end{figure}

In \cite{skuma-morimoto-yokata} Morimoto--Sakuma--Yokota show that
the tunnel number of $K_m$ is super additive, \ie, $t(K_m \# K_m)
= 2t(K_m) + 1$; specifically, $t(K_m) = 1$ and $t(K_m \# K_m) =
3$.  In particular, $K_m$ does not admit a primitive meridian
(Proposition \ref{pro:connecting with a PM}).

We show the following:

\medskip
\noindent {\bf Assertion.} The knots $K_m$ are $(1,2)$-knots;
therefore (since $K_m$ are super additive) by Proposition
\ref{prop:sect 1}, $K_m \# K_m$ admits a primitive meridian.

\medskip
\begin{proof}

$K_m$ does not admit a genus 1, 1-bridge presentation, or it would
contain a primitive meridian (Proposition~\ref{characterization of
knots with PM from bridge position}).  Hence for the proof of the
assertion it suffices to show that $K_m$ admits a genus 1,
2-bridge presentation.

Note that $K(7,17)$ is embedded on a genus 1 unknotted torus in
$S^3$, say $T$.  Let $N(T)$ be a regular neighborhood of $T$,
$N(T) = T \times [0,1]$, where $K(7,17) \subset T \times
\{\frac{1}{2}\}$.  Here we regard the projection onto the second
factor of $T \times [0,1]$ as a height function.  Then we can
isotope $K(7,17) \cup \gamma$ in $N(T)$ (see Figure~8) so
that $K(7,17) = \alpha_1 \cup \beta_1 \cup \alpha_2 \cup \beta_2$,
where the $\alpha_i$'s and $\beta_j$'s are arcs, each $\alpha_i$
and $\beta_j$ share exactly one endpoint, and $\alpha_1 \cap
\alpha_2 = \beta_1 \cap \beta_2 = \emptyset$ ($i,j=1,2$), so that:

\begin{enumerate}
    \item $\alpha_1$ and $\alpha_2$ are monotonic,
    \item $\beta_1$ and $\beta_2$ are vertical, and---
    \item $\gamma$ is embedded in $T \times \{ \frac{1}{2} \}$  and
    bounds a disk $D$ in $T \times \{ \frac{1}{2} \}$ such that
    $D \cap K(7,17)= D \cap (\beta_1 \cup \beta_2)$
consists of two points,
    one in $\beta_1$ and the other in $\beta_2$.
\end{enumerate}

\begin{figure}[ht]
\begin{center}
\includegraphics[width=10cm, clip]{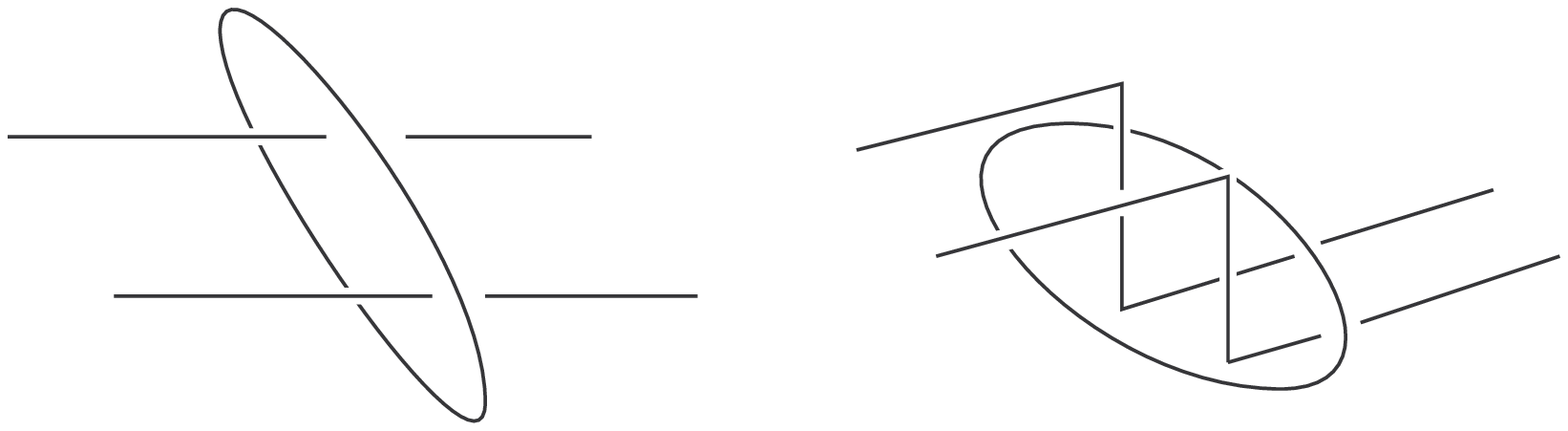}
\end{center}

\begin{center}
Figure 8
\end{center}

\end{figure}

We note that the new position of $K(7,17)$ is a genus 1, 2-bridge
presentation.  Twisting about $\gamma$ changes only the arcs $\beta_1$
and $\beta_2$, but after the twist they remain monotonic; this
completes the proof of the assertion.
\end{proof}

By Theorem \ref{thm:main} and Assertion we have that $gr_t(K_m)
\leq 1/2$.

\providecommand{\bysame}{\leavevmode\hbox
to3em{\hrulefill}\thinspace}
\providecommand{\MR}{\relax\ifhmode\unskip\space\fi MR }
\providecommand{\MRhref}[2]{%
  \href{http://www.ams.org/mathscinet-getitem?mr=#1}{#2}
} \providecommand{\href}[2]{#2}

\end{document}